\RequirePackage{fix-cm}
\documentclass[smallextended]{svjour3}
\smartqed  
\usepackage{graphicx}
\usepackage{epsfig,amsfonts,amssymb,latexsym,amsmath}

 \journalname{my journal}
\begin{document}
\title{A note on optimization in $\mathbb{R}^n$}

\titlerunning{Optimization}        

\author{Fabio Silva Botelho }


\institute{Fabio Silva Botelho \at Department of Mathematics \\
              Federal University of Santa Catarina, SC - Brazil \\
              Tel.: +55-48-3721-3663\\
              \email{fabio.botelho@ufsc.br}  }        


\maketitle

\begin{abstract}
In this article, we develop an algorithm suitable for constrained optimization in $\mathbb{R}^n$. The results are developed through standard tools
of n-dimensional real analysis and basic concepts of optimization. Indeed, the  well known Banach fixed point theorem has a fundamental role in the main result establishment.
  \\
\end{abstract}
\section{Introduction} In this short letter we develop a proximal algorithm for constrained optimization.

Let $f:\mathbb{R}^n \rightarrow \mathbb{R}$ be a $C^2$ class function. Consider the problem of minimizing locally $f$ subject to $g(x) \leq 0,$
where $g:\mathbb{R}^n \rightarrow \mathbb{R}$ is a given $C^2$ class function.

The lagrangian for this problem, denoted by $L:\mathbb{R}^{n+1} \rightarrow \mathbb{R}$ may be expressed by
$$L(x,\lambda)=f(x)+\lambda^2 g(x).$$

We define the proximal formulation for such a problem, denoted by $L_p$ by
$$L_p(x,\lambda,x_k)=f(x)+\lambda^2 g(x)+\frac{K}{2}|x-x_k|^2.$$

\section{The main result}
Linearizing $L_p$, we propose the following procedure for looking for a critical point of such a function:

Consider
\begin{eqnarray}\tilde{L}_{p}(x,\lambda,x_k)&=&f(x_k)+f'(x_k)\cdot (x-x_k)
+\frac{1}{2}[f''(x_k)(x-x_k)]\cdot(x-x_k) \nonumber \\ &&+\lambda^2(g(x_k)+g'(x_k)\cdot(x-x_k))
+\frac{K}{2}|x-x_k|^2.
\nonumber\end{eqnarray}

Hence from $$\frac{\partial \tilde{L}_{p}(x,\lambda,x_k)}{\partial x}=0$$
we obtain,
$$f''(x_k)(x-x_k)+K(x-x_k)+f'(x_k)+\lambda^2g'(x_k)=0,$$
that is,
$$x-x_k=-(f''(x_k)+KI_d)^{-1}(f'(x_k)+\lambda^2g'(x_k)),$$
and therefore $$x(\lambda,x_k)=x_k-(f''(x_k)+KI_d)^{-1}(f'(x_k)+\lambda^2g'(x_k)),$$
where $I_d$ denotes the $n\times n$ identity matrix.

We define $L_1(\lambda,x_k)=\tilde{L}_p(x(\lambda,x_k),x_k,\lambda)$ so that

\begin{eqnarray}L_1(\lambda,x_k)&=&-\frac{1}{2}[(f''(x_k)+KI_d)^{-1}(f'(x_k)+\lambda^2 g'(x_k))] \cdot (f'(x_k)+\lambda^2g'(x_k))
\nonumber \\ &&+f(x_k)+\lambda^2g(x_k) \end{eqnarray}

From $$\frac{\partial L_1(\lambda,x_k)}{\partial \lambda}=0,$$
we get
\begin{equation}\label{t9}[(f''(x_k)+KI_d)^{-1}(f'(x_k)+\lambda^2g'(x_k))]\cdot g'(x_k)\lambda-\lambda g(x_k)=0,\end{equation}
so that we have  two solutions,
$$\lambda_1=0$$ and
\begin{eqnarray}\label{t11}(\lambda_2^1)^2(x_k)=-\left(\frac{[(f''(x_k)+KI_d)^{-1}f'(x_k)] \cdot g'(x_k)-g(x_k)}{[(f''(x_k)+KI_d)^{-1}g'(x_k)] \cdot g'(x_k)}\right).\end{eqnarray}

Observe that if $(\lambda_2^1)^2(x_k)<0$ then $\lambda_2^1(x_k)$ is complex so that, from the condition $\lambda^2\geq 0$, we obtain
 $$\lambda^2(x_k)=\max\{0, (\lambda_2^1)^2(x_k)\}.$$
 Also, from the generalized inverse function theorem  $\lambda^2(x)$ is locally Lipschtzian (see \cite{16,28,12,19} for details). Hence, we may infer that for a given $x_0 \in \mathbb{R}^n$ there exists
 $r>0$ and $\hat{K}_3>0$ such that $$|\lambda^2(x)-\lambda^2(y)|\leq \hat{K}_3|x-y|,$$ $\forall x,y \in B_r(x_0)$.
 With such results in mind, for such an $x_0\in \mathbb{R}^n$,
define $\{x_k\}$ by
$$x_1=x_0-(f''(x_0)+KI_d)^{-1}(f'(x_0)+\lambda^2(x_0)g'(x_0)),$$
$$x_{k+1}=x_k-(f''(x_k)+KI_d)^{-1}(f'(x_k)+\lambda^2(x_k)g'(x_k)),\; \forall k \in \mathbb{N}.$$

Assume   \begin{equation}\label{us18}g(x_0) < 0\end{equation} and there exists $\hat{K}_1$ such that  $|f''(x)| \leq \hat{K}_1, \; \forall x \in B_r(x_0).$

Define $$K_3=\hat{K}_3 \left(\sup_{x \in B_r(x_0)} |g'(x)|\right),$$
$\alpha_1=2K_3/|K-\hat{K}_1|$ and suppose \begin{equation}\label{opt1}f''(x)+\lambda^2(y) g''(x) \geq \alpha_1 (\hat{K}_1+K)I_d\; \forall x,y \in B_r(x_0).\end{equation}

Suppose also $K$ is such that $K>\hat{K}_1$,  $$0<\alpha_1<1,$$
\begin{equation}\label{t19}\left(1-\frac{\alpha_1}{4}\right)I_d \leq ((f''(x)+K I_d)^{-1})(f''(y)+K I_d)\equiv H(x,y) \leq \left(1+\frac{\alpha_1}{4}\right)I_d,
\end{equation} $\forall x,y \in B_r(x_0)$ and \begin{equation}\label{opt8} 0 \leq  \frac{f''(x)+\lambda(y)^2g''(x)}{K-\hat{K}_1} \leq \left(1-\frac{\alpha_1}{2}\right)I_d,
\forall x,y \in B_r(x_0).
\end{equation}
Observe that since $|f''(x)| \leq \hat{K}_1$, we have $$0\leq (K-\hat{K_1})I_d \leq f''(x)+KI_d,$$ so that
\begin{equation}\label{t18}(f''(x)+KI_d)^{-1} \leq \frac{1}{K-\hat{K_1}}I_d,\end{equation}
and
 \begin{equation}\label{opt0}|(f''(x)+KI_d)^{-1}|K_3 \leq \frac{K_3}{|K-\hat{K}_1|}=\frac{\alpha_1}{2},\; \forall x \in B_r(x_0).\end{equation}

Assume $K>0$ is such that $$x_1 \in B_{r(1-\alpha_0)}(x_0)$$ and suppose the induction hypotheses $$x_2,\ldots,x_{k+1} \in B_r(x_0).$$
where $0<\alpha_0<1$ is specified in the next lines.

Note that,
$$x_{k+2}-x_{k+1}=-(f''(x_{k+1})+KI_d)^{-1}(f'(x_{k+1})+\lambda^2(x_{k+1})g'(x_{k+1})),$$
and
$$x_{k+1}-x_{k}=-(f''(x_{k})+KI_d)^{-1}(f'(x_{k})+\lambda^2(x_k) g'(x_{k})),$$
so that,
$$
(f''(x_{k+1})+KI_d)(x_{k+2}-x_{k+1})=-(f'(x_{k+1})+\lambda^2(x_{k+1}) g'(x_{k+1})),$$
and
$$
(f''(x_{k})+KI_d)(x_{k+1}-x_{k})=-(f'(x_{k})+\lambda^2(x_k) g'(x_{k})).$$
Therefore,
\begin{eqnarray}
&&(f''(x_{k+1})+KI_d)(x_{k+2}-x_{k+1})\nonumber \\ &=&(f''(x_{k})+KI_d)(x_{k+1}-x_{k})\nonumber \\ &&
-(f'(x_{k+1})+\lambda^2(x_{k+1}) g'(x_{k+1}))+(f'(x_{k})+\lambda^2(x_k)g'(x_{k})) \nonumber \\ &=&
(f''(x_k)+KI_d)(x_{k+1}-x_k) -(f'(x_{k+1})+\lambda^2(x_{k+1})g'(x_{k+1}))+(f'(x_k)+\lambda^2(x_{k+1})g'(x_k))
\nonumber \\ && -(f'(x_k)+\lambda^2(x_{k+1})g'(x_k))+(f'(x_k)+\lambda^2(x_{k})g'(x_k))
\nonumber \\ &=&(f''(x_{k})+KI_d)(x_{k+1}-x_{k})\nonumber \\ &&
-(f''(\tilde{x}_{k})+\lambda^2(x_{k+1}) g''(\tilde{x}_{k}))(x_{k+1}-x_k)-(\lambda^2(x_{k+1})-\lambda^2(x_k))g'(x_{k})
 \nonumber \end{eqnarray}
where $\tilde{x}_k$ is on the line connecting $x_k$ and $x_{k+1}.$

Thus,
\begin{eqnarray}
x_{k+2}-x_{k+1}&=&(f''(x_{k+1})+KI_d)^{-1}[(f''(x_k)+KI_d)(x_{k+1}-x_k)
\nonumber \\ &&-(f''(\tilde{x}_k)+\lambda^2(x_{k+1})g''(\tilde{x}_k))(x_{k+1}-x_k)
\nonumber \\ &&-(\lambda^2(x_{k+1})-\lambda^2(x_k))g'(x_{k})],
\end{eqnarray}
so that \begin{eqnarray}\label{opt2}
|x_{k+2}-x_{k+1}|&\leq& |H(x_{k+1},x_k)-((f''(x_{k+1})+KI_d)^{-1})(f''(\tilde{x}_k) \nonumber \\ &&+\lambda^2(x_{k+1})g''(\tilde{x}_k))||x_{k+1}-x_k|
\nonumber \\ && +|(f''(x_{k+1})+K\;I_d)^{-1}|K_3|x_{k+1}-x_k|.
\end{eqnarray}

Observe that, from (\ref{opt1}),
$$f''(\tilde{x}_k)+\lambda^2(\tilde{x}_{k+1}) g''(x_k) \geq \alpha_1 (\hat{K}_1+K)I_d \geq \alpha_1(f''(x_{k+1})+K\;I_d),$$
so that $$((f''(x_{k+1})+K\;I_d)^{-1})(f''(\tilde{x}_k)+\lambda^2(x_{k+1}) g''(\tilde{x}_k)) \geq  \alpha_1 I_d.$$

Hence, from this,  (\ref{t19}), (\ref{t18}) and (\ref{opt8}), we obtain
\begin{eqnarray}&&I_d\left(1+\frac{\alpha_1}{4}\right)-\alpha_1\;I_d \nonumber \\ &\geq& H(x_{k+1},x_k)- ((f''(x_{k+1})+K\;I_d)^{-1})(f''(\tilde{x}_k)+\lambda^2(x_{k+1}) g''(\tilde{x}_k)) \nonumber \\ &\geq&
 I_d\left(1-\frac{\alpha_1}{4}\right)- (K\;I_d-\hat{K}_1\;I_d)^{-1})(f''(\tilde{x}_k)+\lambda^2(x_{k+1}) g''(\tilde{x}_k)) \nonumber \\ &\geq&
I_d\left(1-\frac{\alpha_1}{4}\right)-I_d\left(1-\frac{\alpha_1}{2}\right)\nonumber \\ &=& \frac{\alpha_1}{4} I_d \nonumber \\ &\geq& 0,\end{eqnarray}
and therefore,
$$|H(x_{k+1},x_{k})-(f''(x_k)+K\;I_d)^{-1}(f''(\tilde{x}_k)+\lambda^2(x_{k+1}) g''(\tilde{x}_k))| \leq 1-\frac{3 \alpha_1}{4}.$$

On the other hand, from (\ref{opt0}) we have,
$$|(f''(x_k)+K\;I_d)^{-1}| K_3  \leq \frac{\alpha_1}{2}.$$

From (\ref{opt2}) and these last two inequalities, we obtain
$$|x_{k+2}-x_{k+1}| \leq \left(1-\frac{3\alpha_1}{4}+\frac{\alpha_1}{2}\right)|x_{k+1}-x_{k}|=\left(1-\frac{\alpha_1}{4}\right)|x_{k+1}-x_{k}|.$$

Thus, denoting $\alpha_0=1-\alpha_1/4$, we have obtained,
$$|x_{j+2}-x_{j+1}| \leq \alpha_0|x_{j+1}-x_j|, \forall j \in \{1,\cdots,k+1\}$$
so that
\begin{eqnarray}
|x_{j+2}-x_{j+1}| &\leq& \alpha_0|x_{j+1}-x_j|\nonumber \\ &\leq&
\alpha_0^2|x_{j}-x_{j-1}|\nonumber \\ &\leq& \cdots
 \nonumber \\ &\leq&
 \alpha_0^{j+1}|x_1-x_0|,\; \forall j \in \{1,\cdots, k\}.
 \end{eqnarray}
Thus,
\begin{eqnarray}&& |x_{k+2}-x_1|
\nonumber \\ &=& |x_{k+2}-x_{k+1} +x_{k+1}-x_{k}+x_{k}
-x_{k-1}+\cdots + x_2-x_1| \nonumber \\ &\leq&
|x_{k+2}-x_{k+1}| +|x_{k+1}-x_{k}|+\cdots + |x_2-x_1|
\nonumber \\ &\leq& \sum_{j=1}^{k+1} \alpha_0^j|x_1-x_0| \nonumber \\ &\leq&
\sum_{j=1}^{+\infty} \alpha_0^j|x_1-x_0| \nonumber \\ &=&
\frac{\alpha_0}{1-\alpha_0}|x_1-x_0|, \end{eqnarray}
so that
\begin{eqnarray}|x_{k+2}-x_0| &\leq& |x_{k+2}-x_1|+|x_1-x_0| \nonumber \\ &\leq&
\frac{\alpha_0}{1-\alpha_0}|x_1-x_0|+|x_1-x_0| \nonumber \\ &=&
\frac{1}{1-\alpha_0}|x_1-x_0| \nonumber \\ &<&\frac{1}{1-\alpha_0}r(1-\alpha_0) \nonumber \\ &=& r.
\end{eqnarray}

Hence $x_{k+2} \in B_r(x_0),$  and therefore the induction is complete, so that,
$$x_{k}  \in B_r(x_0),\; \forall k \in \mathbb{N}^.$$

Moreover, $\{x_k\}$ is a Cauchy sequence, so that there exists $\tilde{x},$ such that
$$x_k \rightarrow \tilde{x},\; \text{ as } k \rightarrow \infty.$$
Finally
\begin{eqnarray}0&=&\lim_{k \rightarrow \infty} (x_{k+1}-x_k)
\nonumber \\ &=& \lim_{k \rightarrow \infty}[-(f''(x_k)+KI_d)^{-1}(f'(x_k)+\lambda^2(x_k)g'(x_k))]
\nonumber \\ &=&-(f''(\tilde{x})+KI_d)^{-1}(f'(\tilde{x})+\tilde{\lambda}^2g'(\tilde{x})).\end{eqnarray}

Hence, from this and $$det(f''(\tilde{x})+KI_d) \neq 0,$$ we obtain
 $$f'(\tilde{x})+\tilde{\lambda}^2g'(\tilde{x})=0$$
 In such a case, from (\ref{t9}) letting $k \rightarrow \infty,$ we also obtain
 $$\tilde{\lambda}^2g(\tilde{x}) =0.$$

 Thus if $\tilde{\lambda}^2>0$, then $g(\tilde{x})=0.$

 If $\tilde{\lambda}=0,$ then $f'(\tilde{x})=\mathbf{0}$ and $$(\lambda_2^1)(\tilde{x}) \leq 0$$ so that from (\ref{t11}),
 since $(f''(\tilde{x})+KI_d)^{-1}$ is positive definite, letting $k \rightarrow \infty$, we get
 $$g(\tilde{x})=(\lambda_2^1)(\tilde{x})[(f''(\tilde{x})+KI_d)^{-1}g'(\tilde{x})] \cdot g'(\tilde{x}) \leq 0.$$

 That is, in any case, $$g(\tilde{x}) \leq 0.$$

 \begin{remark}\label{t3} For the more general case  with $m_1$ equality  scalar constraints
 $$h_j(x)=0,\forall j \in \{1,\ldots,m_1\}$$ and $m_2$ inequality scalar constraints $$g_l(x) \leq 0,\;\forall l \in \{1,\ldots,m_2\},$$ where
 $h_j,g_l:\mathbb{R}^n \rightarrow \mathbb{R}$ are $C^2$ class functions, $\forall j \in \{1,\ldots,m_1\}$ and
  $\forall l \in \{1,\ldots,m_2\},$ we assume $m_1+m_2<n$ and define the Lagrangian  $L_p$ by
$$L_p(x,\lambda,x_k)=f(x)+\sum_{j=1}^{m_1} (\lambda_h)_j h_j(x) +\sum_{l=1}^{m_2}(\lambda_g)_l^2 g_l(x)+\frac{K}{2}|x-x_k|^2.$$

Linearizing  $L_p$, we propose the following procedure for looking for a critical point of such a function:

Consider
\begin{eqnarray}\tilde{L}_{p}(x,\lambda,x_k)&=&f(x_k)+f'(x_k)\cdot (x-x_k)
+\frac{1}{2}[f''(x_k)(x-x_k)]\cdot(x-x_k) \nonumber \\ &&+\sum_{j=1}^{m_1}(\lambda_h)_j(h_j(x_k)+h_j'(x_k)\cdot(x-x_k))
\nonumber \\ && +\sum_{l=1}^{m_2}(\lambda_g)_l^2(g_l(x_k)+g_l'(x_k)\cdot(x-x_k))
+\frac{K}{2}|x-x_k|^2.
\nonumber\end{eqnarray}

Hence
from $$\frac{\partial \tilde{L}_{p}(x,\lambda,x_k)}{\partial x}=0,$$
we obtain,
$$f''(x_k)(x-x_k)+K(x-x_k)+f'(x_k)+\sum_{j=1}^{m_1} (\lambda_h)_jh_j'(x_k)+\sum_{l=1}^{m_2}(\lambda_g)_l^2g'_l(x_k)=0,$$
that is,
$$x-x_k=-(f''(x_k)+KI_d)^{-1}\left(f'(x_k)+\sum_{j=1}^{m_1} (\lambda_h)_jh_j'(x_k)+\sum_{l=1}^{m_2}(\lambda_g)_l^2g'_l(x_k)\right),$$
and therefore \begin{equation}\label{t1}x(\lambda,x_k)=x_k-(f''(x_k)+KI_d)^{-1}\left(f'(x_k)+\sum_{j=1}^{m_1} (\lambda_h)_jh_j'(x_k)+\sum_{l=1}^{m_2}(\lambda_g)_l^2 g'_l(x_k)\right),\end{equation}
where $I_d$ denotes the $n\times n$ identity matrix.

We define $L_1(\lambda,x_k)=\tilde{L}_p(x(\lambda,x_k),x_k,\lambda),$ so that

\begin{eqnarray}L_1(\lambda,x_k)&=&-\frac{1}{2}\left[(f''(x_k)+KI_d)^{-1}\left(f'(x_k)+\sum_{j=1}^{m_1} (\lambda_h)_jh_j'(x_k)+\sum_{l=1}^{m_2}(\lambda_g)_l^2 g'_l(x_k)\right)\right. \nonumber \\ && \left.\cdot \left(f'(x_k)+\sum_{j=1}^{m_1} (\lambda_h)_jh_j'(x_k)+\sum_{l=1}^{m_2}(\lambda_g)_l^2 g'_l(x_k)\right)\right]
\nonumber \\ &&+f(x_k)+\sum_{j=1}^{m_1}(\lambda_h)_j h_j(x_k)+\sum_{l=1}^{m_2}(\lambda_g)_l^2 g_l(x_k). \end{eqnarray}

From $$\frac{\partial L_1(\lambda,x_k)}{\partial (\lambda_g)_l}=0,$$
we get
\begin{eqnarray}&&\left[(f''(x_k)+KI_d)^{-1}\left(f'(x_k)+\sum_{j=1}^{m_1} (\lambda_h)_jh_j'(x_k)\right.\right.
\nonumber \\ &&\left.\left.+\sum_{l=1}^{m_2}(\lambda^2_g)_l g'_l(x_k)\right)\right]\cdot g_l'(x_k)(\lambda_{g})_l-(\lambda_g)_l g_l(x_k)=0,\end{eqnarray}

From $$\frac{\partial L_1(\lambda,x_k)}{\partial (\lambda_h)_j}=0,$$ we have
\begin{eqnarray}\label{t21}&&\left[(f''(x_k)+KI_d)^{-1}\left(f'(x_k)+\sum_{j=1}^{m_1} (\lambda_h)_jh_j'(x_k)\right.\right.
\nonumber \\ && \left.\left.+\sum_{l=1}^{m_2}(\lambda_g)_l^2 g'_l(x_k)\right)\right]\cdot h_j'(x_k)-h_j(x_k)=0,\end{eqnarray}
$\forall j \in \{1,\ldots,m_1\}.$
Solving the linear system which comprises these last $m_1$ equations and the $m_2$ equations
\begin{eqnarray}\label{t22}&&\left[(f''(x_k)+KI_d)^{-1}\left(f'(x_k)+\sum_{j=1}^{m_1} (\lambda_h)_jh_j'(x_k) \right.\right. \nonumber \\ &&
\left.\left.+\sum_{l=1}^{m_2}(\lambda_g)_l^2 g'_l(x_k)\right)\right]\cdot g_l'(x_k) -g_l(x_k)
= 0,\end{eqnarray} $\forall l \in \{1,\ldots,m_2\},$
 we may obtain a solution $$\left((\lambda_h)_j(x_k),(\lambda_g^1)_l^2(x_k)\right).$$

Thus, to obtain a concerning critical point, we follow the following algorithm.

\begin{enumerate}
\item Choose $x_0 \in \mathbb{R}^n$, $K_{max} \in \mathbb{N}$ ($K_{max}$ is the maximum number of iterations), set $k=0$ and  $e_1 \approx 10^{-5}$.
\item\label{t27} Obtain a solution  $$\left((\lambda_h)_j(x_k),(\lambda_g^1)_l^2(x_k)\right)$$ by solving the linear system
(in $(\lambda_h)_j$ and $(\lambda_g)_l^2$) indicated in (\ref{t21}) and (\ref{t22}).

Observe  that if $(\lambda_g^1)_l^2<0$ then $(\lambda_g^1)_l$ is complex.

  To up-date $\lambda_h$ and $\lambda_g$ proceed as follows:

 \item\label{t7} For each $l \in \{1,\ldots, m_2\}$ if $(\lambda_g)_l^2(x_k)\leq 0$, then  set $(\lambda_g)_l(x_k)=0.$
 \item Define $J=\{l \in \{1,\ldots,m_2\} \text{ such that } (\lambda_g)_l^2(x_k)>0\}.$
 \item Recalculate $(\lambda_h)_j(x_k)$ and the non-zero $(\lambda_g)_l^2(x_k) \text{ for }l \in J$ through the solution of the linear system (in $(\lambda_h)_j$ and $(\lambda_g)_l^2$)
 \begin{eqnarray}&&\left[(f''(x_k)+KI_d)^{-1}\left(f'(x_k)+\sum_{j=1}^{m_1} (\lambda_h)_jh_j'(x_k)\right.\right.
\nonumber \\ && \left.\left.+\sum_{l \in J}(\lambda_g)_l^2 g'_l(x_k)\right)\right]\cdot h_j'(x_k)-h_j(x_k)=0,\end{eqnarray}
$\forall j \in \{1,\ldots,m_1\}$
and
\begin{eqnarray}&&\left[(f''(x_k)+KI_d)^{-1}\left(f'(x_k)+\sum_{j=1}^{m_1} (\lambda_h)_jh_j'(x_k)\right.\right. \nonumber \\ && \left.\left.+\sum_{l \in J}(\lambda_g)_l^2 g'_l(x_k)\right)\right]\cdot g_l'(x_k)-g_l(x_k)
=0,\end{eqnarray} $\forall l \in J.$
\item If $(\lambda_g)_l^2(x_k) \geq 0,\; \forall l \in \{1,\ldots,m_2\}$, then  go to \ref{t23}, otherwise go to item \ref{t7}.
\item\label{t23} Up-date $x_k$ through the equation
  \begin{eqnarray}\label{t1}x_{k+1}&=&x_k-(f''(x_k)+KI_d)^{-1}\left(f'(x_k)+\sum_{j=1}^{m_1} (\lambda_h)_j(x_k)h_j'(x_k)
\right. \nonumber \\ && \left.+\sum_{l=1}^{m_2}(\lambda_g)_l^2(x_k) g'_l(x_k)\right).\end{eqnarray}
 \item If $|x_{k+1}-x_k| < e_1$ or $k> K_{max}$, then stop, otherwise $k:= k+1$ and go to  \ref{t27}.
 \end{enumerate}
\end{remark}
\section{Conclusion} In this article we have developed an algorithm for constrained optimization in $\mathbb{R}^n$. We prove the main result only for
the special case of a single scalar inequality constraint. However, we highlight the proof of a more general result involving equality and inequality
constraints may be developed in a similar fashion, as indicated in remark \ref{t3}. We postpone the presentation of the formal details for such a more
general case for a future work.


\end{document}